\documentclass[12pt]{article}

\def\Ker{{\rm Ker}}
\def\Im{{\rm Im}}
\def\Fix{{\rm  Fix}}
\def\Fixe{{\rm Fix}\,\vp}
\def\vp{{\varepsilon}}
\def\w{\wedge }
\def\we{\,{\wedge}_\varepsilon\,}
\def\ve{\,{\vee}_\varepsilon\,}
\def\weu{\,{\wedge}_{\varepsilon_1}\,}
\def\veu{\,{\vee}_{\varepsilon_1}\,}
\def\wed{\,{\wedge}_{\varepsilon_2}\,}
\def\ved{\,{\vee}_{\varepsilon_2}\,}
\def\cale{{\cal E}}

\def\bit{\begin{itemize}}
\def\eit{\end{itemize}}
\def\barr{\begin{array}}
\def\earr{\end{array}}
\def\calel{$\cale$-lattice}
\def\fwg{following}
\def\vspp{\vspace*{1,6mm}\\ }
\def\ord{{\rm ord}}

\def\Z{{\rlap{$\kern2pt{\rm Z}$}{\rm Z}\,}}
\def\bld#1#2{{\buildrel{#1}\over{#2}}}
\def\st#1#2{{\mathrel{\mathop{#2}\limits_{#1}}{}\!}}
\def\stb#1#2#3{{\st{{#1}}{\bld{{#2}}{#3}}{}\!}}
\def\xmare#1#2{\stb{#1}{#2}{\mbox{\Huge$\times$}}}
\def\a{{\alpha}}
\def\b{{\beta}}
\def\g{{\gamma}}
\def\nn{{\rm I\!N}}
\def\nns{{\rm I\!N}^*}
\def\dd{\displaystyle}

\title{{\bf An $\cale$-lattice structure associated\\ to some classes of finite groups}\footnote{This work has been supported by the
research grant GAR 88/2007-2008}}
\author{Marius T\u arn\u auceanu\\
{\small Faculty of  Mathematics, "Al.I. Cuza" University,
Ia\c si, Romania}\\
{\small e-mail: tarnauc@uaic.ro}}
\date{October 1, 2008}

\begin{document}
\maketitle

\begin{abstract} In the present paper we introduce and study  a canonical $\cale$-lattice structure on the set of element  orders of some finite groups. We show that a finite abelian group is uniquely determined by this canonical $\cale$-lattice. \smallskip
\end{abstract}

\section{Preliminaries}

The notion of $\cale$-lattice has been introduced in \cite{6}. So,
given a nonvoid set $L$ and a map $\vp:L\to L$, we denote by
$\Ker\,\vp$ the kernel of $\vp$ (i.e. $\Ker\,\vp=\{(a,b)\in L\times
L\mid \vp(a)=\vp(b)\}$, by $\Im\,\vp$ the image of $\vp$  (i.e.
$\Im\,\vp=\{\vp(a)\mid a\in L\})$ and by $\Fixe$ the set consisting
of all fixed points of $\vp$ (i.e. $\Fixe=\{a\in L\mid\vp(a)=a\})$.
We say that $L$ is an {\it $\cale$-lattice} (relative to $\vp)$ if
there exist two binary operations $\we,\ve$ on $L$ which satisfy the
\fwg\ properties: \bit\item[a)] $a\we(b\we c){=}(a\we b)\we c,\
a\ve(b\ve c){=}(a\ve b)\ve c,$ for all $a,b,c\in L;$
\item[b)] $a\we b=b\we a,\ a\ve b=b\ve a,$ for all $a,b\in L;$
\item[c)] $a\we a=a\ve a=\vp(a),$ for any $a\in L;$
\item[d)] $a\we (a\ve b)=a\ve(a\we b)=\vp(a),$ for all $a,b\in L.$\eit

Clearly, in an \calel\ $L$ (relative to $\vp)$ the map $\vp$  is
idempotent and $\Im\,\vp=\Fixe.$ Moreover, for any $a,b\in L$, we
have:
$$\barr{l}
a\we\vp(a)=a\ve\vp(a)=\vp(a),\vspp a\we\vp(b)=\vp(a)\we
b=\vp(a)\we\vp(b)=\vp(a\we b),\vspp a\ve\vp(b)=\vp(a)\ve
b=\vp(a)\ve\vp(b)=\vp(a\ve b).\earr$$ Also, note that the set
$\Fixe$ is closed under the binary  operations $\we,\ve$ and,
denoting by $\w,\vee$ the restrictions of $\we,\ve$ to $\Fixe$, we
have that $(\Fixe,\w,\vee)$ is a lattice. The connection between the
\calel\ concept and the lattice concept is very powerful. Thus, if
$(L,\we,\ve)$ is an \calel\ and $\sim$ is an equivalence relation on
$L$ such that $\sim\,\subseteq\Ker\,\vp$, then the factor set
$L/\!\sim$ is a lattice isomorphic to the lattice $\Fixe$.
Conversely, if $L$ is a nonvoid set and $\sim$ is an equivalence
relation on $L$ having the property that the factor set $L/\!\sim$ is a lattice, then the set $L$ can be endowed with an \calel\ structure
(relative to a map $\vp:L\to L)$ such that
$\sim\,\subseteq\Ker\,\vp$ and $L/\!\sim\,\cong\Fixe.$

If $(L,\we,\ve)$ is an \calel\ and for every $x\in L$ we denote  by
$[x]$ the equivalence class of $x$ modulo $\Ker\,\vp$ (i.e.
$[x]=\{y\in L\mid\vp(x)=\vp(y)\}),$ then we have $a\we
b\in[\vp(a)\w\vp(b)]$ and $a\ve b\in[\vp(a)\vee\vp(b)]$, for all
$a,b\in L$. We say that $L$ is {\it canonical \calel} if $a\we b,$
$a\vee b\in \Fixe,$ for all $a,b\in L$. Three fundamental  types of
canonical \calel s have been identified and studied in \cite{6}.

Let $(L_1,\weu,\veu),(L_2,\wed,\ved)$ be  two \calel s and,  for
every element $a\in L_i,$ denote by $[a]_i$ the equivalence class of
$a$~modulo $\Ker\,\vp_i$, $i=1,2.$ A map $f:L_1\to L_2$ is called an
{\it\calel\ homomorphism} if: \bit\item[a)] $f\circ\vp_1=\vp_2\circ
f;$
\item[b)] for all $a,b\in L_1,$ we have:
\bit\item[i)] $f(a\weu b)=f(a)\wed f(b);$
\item[ii)] $f(a\veu b)=f(a)\ved f(b).$\eit\eit
Moreover, if the map $f$ is one-to-one and onto, then it is called
an {\it\calel\ isomorphism}. \calel\ isomorphisms between canonical
\calel s have been investigated in \cite{7}. Mention here only the
\fwg\ characterization of these, established in Proposition 1,
\S~2.1: a map $f:L_1\to L_2$ is an \calel\ isomorphism if and only
if its restriction to the set $\Fix\,\vp_1$ is a lattice isomorphism
from $\Fix\,\vp_1$ to $\Fix\,\vp_2$ and $f|[a]_1:[a]_1\to [f(a)]_2$
is one-to-one and onto, for  each $a\in\Fix\,\vp_1.$

An interesting example of a canonical \calel\ structure  on the set
of subgroups of a group $G$ has been presented in \cite{7}. Another
canonical \calel\ can be associated to some classes of finite groups
$G$, too. Its study is the main goal of this paper.

Most of our notations are standard and will usually not be  repeated
here. Basic definitions and results on lattices can be found in
\cite{1} or \cite{2}. For group theory concepts we refer the reader
to \cite{3} and \cite{4}.

\section{Main results}

Let $G$ be a finite group of order $n$ and $L_n$ be the lattice
consisting of all divisors of $n$. We define on $G$ the \fwg\
equivalence relation
$$a\sim a'\mbox{\ \ iff\ \ } \ord(a)=\ord(a')$$
and take a set of representatives $\{a_i\mid i\in I\}$  for the
equivalence classes modulo~$\sim$. There exist many classes of
finite groups such that the set $\{\ord(a_i)\mid i\in I\}$ forms a
sublattice of $L_n$, as finite  $p$-groups or finite abelian groups.
For such a  group $G$, this induces a lattice structure on
$\{a_i\mid i\in I\}$, which is isomorphic to that of
$\{\ord(a_i)\mid i\in I\}$. As we have seen in Section 1, $G$
becomes a canonical \calel\ with respect to the idempotent map
$$\vp:G\to G,\ \vp(a)=a_i\mbox{\ \ iff\ \ }a\in[a_i],$$where the binary operations $\we,\ve$ are defined by
$$a\we b=\vp(a)\w\vp(b)\mbox{\ \ and\ \ }a\ve b=\vp(a)\vee\vp(b).$$
We shall call it  the {\it order canonical \calel} associated to $G$.\bigskip

Our first aim is to  give a detailed description of the above
\calel\  in the case when $G$ is a finite abelian group. Then, from
the fundamental theorem of finitely generated abelian groups, there
are (uniquely determined by $G$) the numbers $m\in\nns$,
$d_1,d_2,...,d_m\in\nn\setminus\{0,1\}$ satisfying
$d_1/d_2/.../d_m$, $\dd\prod_{i=1}^md_i=n$ and
$G\cong\xmare{i=1}m\Z_{d_i}.$ Clearly, the element orders of $G$
are all divisors of the exponent $d_m$ of $G$, therefore the \fwg\
lattice isomorphism holds:
$$\Fixe\cong L_{d_m}.\leqno(*)$$
In order to complete our description we need to count the elements
of an arbitrary equivalence  class $[a_i]$ modulo~$\sim$, i.e. to
determine the number of elements of a given order in $G$. First of
all, we focus on the particular situation in which $G$ is a finite
abelian $p$-group. Then it has a direct decomposition of type
$G\cong\xmare{i=1}m\Z_{p^{\a_i}},$ where $p$ is a prime and
$1\le\a_1\le\a_2\le...\le\a_m$. A recurrence method to compute the
number of elements of order $p^\a$ $(\a\in\nn)$ in $G$ is indicated
in the next lemma.

\bigskip\noindent{\bf Lemma 1.} {\it Let $p$ be  a prime,
$(\a_m)_{m\in\nns}$
be a chain of positive integers such that
$1\le\a_1\le\a_2\le...\le\a_m\le...$  and $(f^p_m)_{m\in\nns}$ be
the chain of functions defined by:
$$\barr{l}
f^p_m:\nn\to\nn,\vspp
f^p_m(\a)=\left|\left\{x\in\xmare{i=1}m\Z_{p^{\a_i}}\mid\ord(x)=p^\a\right\}\right|,\ \ (\forall)\a\in\nn.\earr$$
Then $f^p_m(0)=1,\ f^p_m(\a)=0$ for all $\a>\a_m$ and the chain $(f^p_m)_{m\in\nns}$ satisfies the equality
$$f^p_{m+1}(\a)=p^\a g^p_m(\a)-p^{\a-1}g^p_{m}(\a-1),\ \ (\forall)\a\in\nns,$$where
$$g^p_m(\a)=\left\{\barr{ll}
p^{m \a},&0\le\a\le\a_1\vspp
p^{(m-1)\a+\a_1},&\a_1\le\a\le\a_2\\\vdots\\
p^{\a_1+\a_2+...+\a_m},&\a_m\le\a.\earr\right.$$}

\noindent{\bf Proof.} Let
$x=(x_1,x_2,...,x_{m+1})\in\xmare{i=1}{m+1}\Z_{p^{\a_i}}.$  Then
$\ord(x)$ is the maximum  between the order of $(x_1,x_2,...,x_m)$
in $\xmare{i=1}m\Z_{p^{\a_i}}$ and the order of $x_{m+1}$ in
$\Z_{p^{\a_{m+1}}},$ and so, $\ord(x)=p^\a$ iff:
$$ \left\{\barr{ll}
\ord(x_1,x_2,...,x_m)=p^\a\mbox{ and }\ord(x_{m+1})=p^\b\mbox{ for some $\b\le\a$}\vspp
\mbox{or}\vspp
\ord(x_1,x_2,...,x_m)=p^\b\mbox{ for some $\b<\a$ and  $\ord(x_{m+1})=p^\a.$}\earr\right.$$
This shows that the \fwg\  recurrence relation holds
$$\barr{l }
f^p_{m+1}(\a) {=} p^\a f^p_m(\a){+}(p^{\a}{-}p^{\a-1})(f^p_m(0){+}f^p_m(1){+}...{+}f^p_m(\a{-}1))=\vspp
 {=} p^\a(f^p_m(0){+}f^p_m(1){+}...{+}f^p_m(\a)){-}p^{\a-1}(f^p_m(0){+}f^p_m(1){+}...{+}f^p_m(\a{-}1)),\earr\leqno(1)$$
 for all $\a\in\nns$. Putting $g^p_m(\a)=\dd\sum_{\b=0}^\a f^p_m(\b)$, the equality (1) becomes:
$$f^p_{m+1}(\a)=p^\a g^p_m(\a)-p^{\a-1}g^p_m(\a-1),\ \ (\forall)\a\in\nns.\leqno(2)$$
Summing up these equalities for $\a=1,2,...,$ we obtain:
$$g^p_{m+1}(\a)=p^\a g^p_m(\a),\ \ (\forall)\a\in\nn.\leqno(3)$$Since
$$\barr{lcl}
g^p_1(\a)&=&\dd\sum_{\b=0}^\a f^p_1(\b)=1+\dd\sum_{\b=1}^{\min\{\a,\a_1\}}f^p_1(\b)=1+\dd\sum_{\b=1}^{\min\{\a,\a_1\}}(p^\b-p^{\b-1})=\vspp
&=&\left\{\barr{ll}
p^\a,&\ 0\le\a\le\a_1\vspp p^{\a_1},&\ \a_1\le\a,\earr\right.\earr$$
by (3) we infer that $g^p_m(\a)$ is of the indicated form.\hfill\rule{1,5mm}{1,5mm}

\bigskip\noindent{\bf Example.} We want to determine the number of
 elements of a given order in the group $\Z_4\times\Z_{16}.$ By applying Lemma 1 for $m=1,$ we obtain:
$$f^p_2(\a)=\left\{\barr{ll}
1,&\a=0\vspp
p^{2\a}-p^{2(\a-1)},&1\le\a\le\a_1\vspp
p^{\a+\a_1}-p^{\a-1+\a_1},&\a_1<\a\le\a_2\vspp
0,&\a_2<\a.\earr\right.$$
Taking $p=2,\ \a_1=2,\ \a_2=4,$ it results that
$$f^2_2(\a)=\left\{\barr{ll}
1,&\a=0\vspp 2^{2\a}-2^{2(\a-1)},&1\le\a\le2\vspp 2^{\a+2}-2^{\a+
1},&2<\a\le4\vspp 0,&4<\a,\earr\right.$$ and thus
$\Z_4\times\Z_{16}$ has 1 element of order $2^0,$ $2^2-2^0=3$
elements of order $2^1$, $2^4-2^2=12$ elements of order $2^2$,
$2^5-2^4=16$ elements of order $2^3$ and $2^6-2^5=32$ elements of
order $2^4.$\bigskip

Next, we return to the  general case in which $G$ is a finite
abelian group. Because such a group can be uniquely written as a
direct product of finite abelian $p$-groups, by using Lemma 1 we are
able to establish an explicit formula for the number of elements of
an arbitrary order in $G$.

\bigskip\noindent{\bf Theorem 2.} {\it Let $G$ be a finite
abelian group of order $n$, $\{p_1,p_2,...,p_k\}$  be the set of all
distinct prime divisors of $n$ and $G\cong\xmare{i=1}kG_i$ be the
direct decomposition of $G$ as product of finite abelian $p$-groups
$(G_i=p_i$-group, $i=\overline{1,k})$. Then, for any
$(\b_1,\b_2,...,\b_k)\in\nn^k$, the number of elements of order
$p_1^{\b_1}p_2^{\b_2}...p_k^{\b_k}$ in $G$ is
$$\prod_{i=1}^kf^{p_i}_{m_i}(\b_i),$$where $m_i$ denotes  the number of direct factors of $G_i$, $i=\overline{1,k}$.}

\bigskip\noindent{\bf Proof.} Since every element $x$ of $G$ can be
uniquely written as a product $x_1x_2...x_k$ with $x_i\in G_i$,
$i=\overline{1,k}$, and $(\ord(x_i),\ord(x_j))=1$ for all $i\ne j
 =\overline{1,k}$,  we get $\ord(x){=}\dd\prod_{i=1}^k\ord(x_i)$.
 So, $\ord(x){=}p_1^{\b_1}p_2^{\b_2}...p_k^{\b_k}$ if and only if
 $\ord(x_i){=}p_i^{\b_i},$ $i{=}\overline{1,k}$. Now, our statement
 follows immediately from Lemma~1.\hfill\rule{1,5mm}{1,5mm}

\bigskip\noindent{\bf Example.} We want to determine the number of
elements of  order 120 in the group $\Z_{12}\times\Z_{720}.$ First
we write $G$ as a product of finite abelian $p$-groups
$$G=G_1\times G_2\times G_3,$$
where $G_1=\Z_4\times\Z_{16},$ $G_2=\Z_3\times\Z_{9}$ and
$G_3=\Z_5.$ Since $120=2^33^15^1$, we must compute the numbers
$f^2_2(3),f^3_2(1),f^5_1(1)$  of elements of orders $2^3,3^1,5^1$ in
$G_1,G_2,G_3$, respectively. We already know that $f^2_2(3)=16.$
Similarly, by Lemma 1, we find $f^3_2(1)=8$ and $f^5_1(1)=4$. Hence
$\Z_{12}\times\Z_{720}$ has $16\cdot8\cdot4=512$ elements of order
120.\bigskip

From Theorem 2 we easily obtain the number of elements whose orders are prime powers in a finite abelian group.

\bigskip\noindent{\bf Corollary 3.} {\it Under the hypotheses of
Theorem 2, for each $i\in\{1,2,...,k\}$, the number of elements of
order $p_i^\a$ in $G$ is $f^{p_i}_{m_i}(\a)$. In particular, $G$
possesses $f^{p_i}_{m_i}(1)=p_i^{m_i}-1$ elements of order
$p_i$.}\bigskip

Since the structure of the order canonical $\cale$-lattice of a
finite abelian group is completely determined, it appears the
following question: what can be said about two finite abelian groups
whose order canonical \calel s are isomorphic? In order to answer
this question, we need first to prove an auxiliary result.

\bigskip\noindent{\bf Lemma 4.} {\it Let $G$ and $G'$ be
two finite abelian $p$-groups. If the numbers of elements of order
$p^\a$ in $G$ and $G'$ are equal for all $\a\in\nn$, then $G$ and
$G'$ are isomorphic.}

\bigskip\noindent{\bf Proof.} Suppose that
$G=\xmare{i=1}{r+1}\Z_{p^{\a_i}}$ and
$G'\cong\xmare{i=1}{s+1}\Z_{p^{\b_i}}$, where
$1\le\a_1\le\a_2\le...\le\a_{r+1}$  and
$1\le\b_1\le\b_2\le...\le\b_{s+1}.$ By Corollary 3, $G$ and $G'$
have $p^{r+1}-1$ and $p^{s+1}-1$ elements of order $p$,
respectively,  so $r+1=s+1.$ Suppose that $\a_1\ne\b_1$ and assume
$\a_1<\b_1$. Then, by Lemma 1, we obtain that the numbers of
elements of order $p^{\a_1+1}$ in $G$ and $G'$ are
$p^{(r+1)\a_1+r}-p^{(r+1)\a_1}$ and
$p^{(r+1)(\a_1+1)}-p^{(r+1)\a_1},$ fact which contradicts our
hypothesis. Thus $\a_1=\b_1$. Now, it is clear that a standard induction
argument will show that $\a_i=\b_i$, for all $i=\overline{1,r+1}$.
Hence $G\cong G'$.\hfill\rule{1,5mm}{1,5mm}\bigskip

Finally, we are able to show that the order canonical \calel\  of a
finite abelian group determines uniquely the structure of this
group, and so to answer the above question.

\bigskip\noindent{\bf Theorem 5.} {\it Two finite abelian groups $G$
and $G'$ are isomorphic if and only if their order canonical \calel
s are isomorphic.}

\bigskip\noindent{\bf Proof.} Clearly, the order canonical \calel s
of two isomorphic  finite abelian groups are also isomorphic.
Conversely, assume that the order canonical \calel s of
$G\cong\xmare{i=1}m\Z_{d_i}$ and $G'\cong\xmare{i=1}{m'}\Z_{d'_i}$
(where $d_i,$ $i=\overline{1,m},$ and $d'_i,$ $i=\overline{1,m'}$,
are given by the fundamental theorem of finitely generated abelian
groups) are isomorphic, and take an isomorphism $f$ between these.
Then, by using the lattice isomorphism $(*)$, we get:
$$L_{d_m}\cong L_{d'_{m'}}.\leqno(4)$$
This shows that $n=|G|$ and $n'=|G'|$ have the same number of prime
divisors. Put $n=p_1^{\mu_1}p_2^{\mu_2}...p_k^{\mu_k}$ and
$q_1^{\eta_1}q_2^{\eta_2}...q_k^{\eta_k}$, where all $p_i$ and $q_i$
are prime, $i=\overline{1,k}.$ The relation (4) implies also that
$d_m=p_1^{\g_1}p_2^{\g_2}...p_k^{\g_k}$ and
$d'_{m'}=q_1^{\g_1}q_2^{\g_2}...q_k^{\g_k}$ for some
$\g_1,\g_2,...,\g_k\in\nns$. Then $f(p_i^\a)=q^\a_i$, for all
$i=\overline{1,k}$ and $\a\in\nns.$ Because an \calel\ isomorphism
induces a bijection between the corresponding equivalence classes
modulo~$\sim$, it results that the numbers of elements of orders
$p_i^\a$ and $q_i^\a$ in $G$ and $G'$, respectively, are equal.

On the other hand, $G$ and $G'$ can be written as  direct  products
of finite abelian $p$-groups
$$G\cong\xmare{i=1}kG_i,\ \ G'\cong\xmare{i=1}kG'_i,$$
where $|G_i|=p_i^{\mu_i}$ and $|G'_i|=q_i^{\eta_i}$,
$i=\overline{1,k}.$ From Corollary 3 we know that the numbers of
elements of orders $p_i$ in $G_i$ and $q_i$ in $G'_i$ are
$p_i^{m_i}-1$ and $q_i^{m'_i}-1$ $(m_i$ and $m'_i$ being the numbers
of direct factors in $G_i$ and $G'_i)$. It follows that
$p_i^{m_i}-1=q_i^{m'_i}-1$, so $p_i=q_i$, $i=\overline{1,k}.$ Now,
Lemma 4 implies that $G_i\cong G'_i$, for all $i=\overline{1,k},$
and hence $G\cong G'.$\hfill\rule{1,5mm}{1,5mm}\bigskip

We finish our paper by  the remark that the detailed study of order
canonical \calel s can be extended from finite abelian groups to
another classes of finite groups.

\bigskip\noindent{\bf Open problem.} {\it Find other finite groups for
which the previous canonical \calel\ structure can be defined
$($i.e. the orders of elements form, with respect to divisibility, a
lattice$)$ and described $($i.e. we can determine the number of
elements of an arbitrary given order$)$.}


\begin{thebibliography}{7}
\bibitem{1}  Birkhoff, G., {\it Lattice theory}, Amer. Math. Soc., Providence, R.I., 1967.
\bibitem{2}  Gr\"atzer, G., {\it General lattice theory}, Academic Press, New York, 1978.
\bibitem{3} Suzuki, M., {\it Group theory}, I, II, Springer-Verlag, Berlin, 1982, 1986.
\bibitem{4}  \c{S}tef\u anescu, M.,  {\it Introduction to group theory} (Romanian), Ed. Univ. "Al.I. Cuza", Ia\c si, 1993.
\bibitem{5} T\u arn\u auceanu, M.,  {\it Groups determined by posets of subgroups}, Ed. Matrix Rom, Bucure\c sti, 2006.
\bibitem{6}  T\u arn\u auceanu, M.,  {\it $\cale$-lattices}, Italian Journal of Pure and Applied Ma\-the\-ma\-tics 22 (2007), 27-38.
\bibitem{7}  T\u arn\u auceanu, M.,  {\it On isomorphisms of canonical $\cale$-lattice}, Fixed Point Theory 8 (2007), no. 1, 131-139.
\end{thebibliography}
\end{document}